\documentclass[11pt,twoside]{amsart}
\usepackage{amsmath, amsthm, amscd, amsfonts, amssymb, float, graphicx, 
, soul }
\usepackage[bookmarksnumbered, 
, plainpages]{hyperref}
\usepackage{graphicx}
\usepackage{epstopdf}
\usepackage{rotating}
\input xy
\xyoption{all}
\sethlcolor{green}
\newtheorem{thm}{Theorem}[section]

\newtheorem{lem}[thm]{Lemma}
\newtheorem{prop}[thm]{Proposition}

\newtheorem{remark}{Remark}

\newtheorem{ex}[thm]{Example}

\numberwithin{equation}{section}

\def\Cay{{\mathrm {Cay}}}
\begin{document}
\openup 1.9\jot

\vspace{1.3 cm}
\title{Some properties of Cayley signed graphs on finite abelian groups}
\author{Mohammad~A.~Iranmanesh$^*$ and Nasrin Moghaddami}

\thanks{{\scriptsize
\hskip -0.4 true cm MSC(2010): Primary: 20E99.
\newline Keywords: Abelian group, balancing, Cayley graph, clustering, signed graph, sign-compatible.\\
$*$Corresponding author}}
\maketitle
\begin{abstract}
Let $\Sigma=(\Gamma, \sigma)$ is a signed graph(or sigraph in short), where $\Gamma$ is a 
underlying graph of $\Sigma$ and $\sigma:E\longrightarrow \{+, -\}$ is a function. Consider $\Gamma=\Cay(\mathbb{Z}_{p_{1}}\times \mathbb{Z}_{p_{1}^{\alpha_{1}}p_{2}^{\alpha_{2}} \ldots p_{k}^{\alpha_{k}}}, \Phi)$, where all $p_{1}, p_{2}, \ldots, p_{k}$ are distinct prime factors and $\Phi=\varphi_{p_{1}}\times\varphi_{p_{1}^{\alpha_{1}}p_{2}^{\alpha_{2}} \ldots p_{k}^{\alpha_{k}}}$. For any positive
integer $n$, $\varphi_{n}=\{\ell| 1\leq \ell<n, \gcd(\ell, n)=1\}$. Motivated by \cite{s14}, 
we will investigate balancing in $\Sigma$ and $L(\Sigma)$, clusterability and sign-compatibility of $\Sigma$.


\end{abstract}
\maketitle
\section{Introduction}
\label{intro}
Let $(G, \cdot)$ be a group and $S=S^{-1}$ be a non empty subset of $G$ not containing the
identity element e of $G$. The Cayley graph $\Cay(G, S)$ is the simple graph having vertex set $G$ and edge set
$\{\{v, vs\} | v\in G, s \in S \}$. To find enough information about Cayley graphs, refer to 
books Biggs \cite{s12} and Godsil, Royle \cite{godsil}.

A signed graph(or sigraph in short) is an ordered
pair $\Sigma=(\Gamma, \sigma)$ where $\Gamma=(V, E)$ is a graph called the underlying graph of $\Sigma$
and $\sigma:E\longrightarrow \{+, -\}$ is a function. $\Sigma$ is all-positive (all-negative) if all its edges are positive
(negative). Moreover, it is said to be homogeneous if it is either all-positive or all-negative
and heterogeneous otherwise. $d^{-}(v)$ ($d^{+}(v)$) represents the number of negative (positive) edges incident at $v$ in $\Sigma$. 
A marked sigraph is an ordered pair $\Sigma_{\mu}= (\Sigma, \mu)$, where $\Sigma=(\Gamma, \sigma)$ is a sigraph
and $\mu: V(\Sigma)\longrightarrow \{+, -\}$ is a function, called a marking of $\Sigma$ (see \cite{s5, s6, s7}).

Behzad and Chartrand \cite{s10} defined line sigraph
$L(\Sigma)$ as the sigraph in which the edges of $\Sigma$ are represented as vertices, two of these
vertices are defined adjacent whenever the corresponding edges in $\Sigma$ have a vertex in
common, any such edge ef is defined to be negative whenever both $e$ and $f$ are
negative edges in $\Sigma$. A positive cycle in $\Sigma$ is a cycle which contains even number of negative edges. 
A sigraph $\Sigma$ is a balanced sigraph if every cycle in $\Sigma$ is positive.
The first classification of balanced graphs was done by Harary in 1953 \cite{s8}.

A clusterable sigraph is a sigraph in which its vertex set is partitioned into pairwise disjoint subsets in 
such a way that whenever there are positive edges they lie in the same subset and negative edges lie 
in the different subsets. These subsets are called clusters (see \cite{s2}). A sigraph is called 
sign-compatible \cite{s3} if it is possible to mark the vertices of the sigraph in such a manner 
that for every negative edge there are negatively marked vertices on its both ends and for any
 positive edge its ends are not assigned negative, sign-incompatible otherwise.

Let $\Sigma=(\Gamma, \sigma)$ is a sigraph whose underlying graph is  
$\Gamma=\Cay(\mathbb{Z}_{p_{1}}\times \mathbb{Z}_{p_{1}^{\alpha_{1}}p_{2}^{\alpha_{2}} \ldots p_{k}^{\alpha_{k}}}, \Phi)$ 
where $p_{1}, p_{2}, \ldots, p_{k}$ are distinct prime factors and $\Phi=\varphi_{p_{1}}\times\varphi_{p_{1}^{\alpha_{1}}p_{2}^{\alpha_{2}} \ldots p_{k}^{\alpha_{k}}}$ and also for an edge $ab$ of $\Sigma$, the function $\Sigma$ is defined as follows \cite{s14}.
 \begin{equation*}
\sigma(ab) = \left\lbrace \begin{array}{ll}
+  & ~~~~~ if ~~a\in \Phi ~~or~~ b\in\Phi, \vspace{0.5cm}\\ 
-  &~~~~~ otherwise.
\end{array} \right.
 \end{equation*}
By $\mathbb{Z}_{n}$ we denote the abelian group of order $n$.
For $p_{1}=2$, the signed graph $\Sigma$, is a disconnected sigraph with two connected
components, say $\Sigma_{1}$ and $\Sigma_{2}$, where
$V(\Sigma_{1})=\{(1, v)| v~is~odd\}\cup \{(0, v)| v~is~even\}$ and
$V(\Sigma_{2})=\{(0, v)| v~is~odd\}\cup \{(1, v)| v~is~even\}$.
Since every Cayley graph $\Cay(G, S)$ is $|S|$-regular (see for example \cite{godsil}),
we find that $\Sigma$ is $|\Phi|$-regular.

In this article some properties of the signed graph $\Sigma$, such as balancing, clusterability 
and sign-compatibility, are investigated.

\section{balance in sigraph $\Sigma$ and line sigraph $\Sigma$}\label{sec:2}
In this chapter, we examined the number of positive and negative edges for 
$\Sigma=(\Gamma, \sigma)$, where 
$\Gamma=\Cay(\mathbb{Z}_{p_{1}}\times \mathbb{Z}_
{p_{1}^{\alpha_{1}}p_{2}^{\alpha_{2}} \ldots p_{k}^{\alpha_{k}}}, \Phi)$ and 
one of the prime factors say $p_{i}$ is $2$ and $\alpha_{j}\geq 1$ for any $j=1,2, \ldots, k$. 
Moreover we find the number of positive and negative edges of 
$\Gamma_{1}=\Cay(\mathbb{Z}_{p}\times \mathbb{Z}_{p^{\alpha}}, \Phi)$, where $\alpha\geq 1$, $p\geq 3$,
 and $\Gamma_{2}=\Cay(\mathbb{Z}_{p}\times \mathbb{Z}_{pq}, \Phi)$, where $p, q\geq 3$ 
 are prime numbers. Next we study the
  balanced property of $\Sigma$ and $L(\Sigma)$.

\begin{figure}[h]
\begin{center}
\includegraphics[width=14cm , height=6cm]{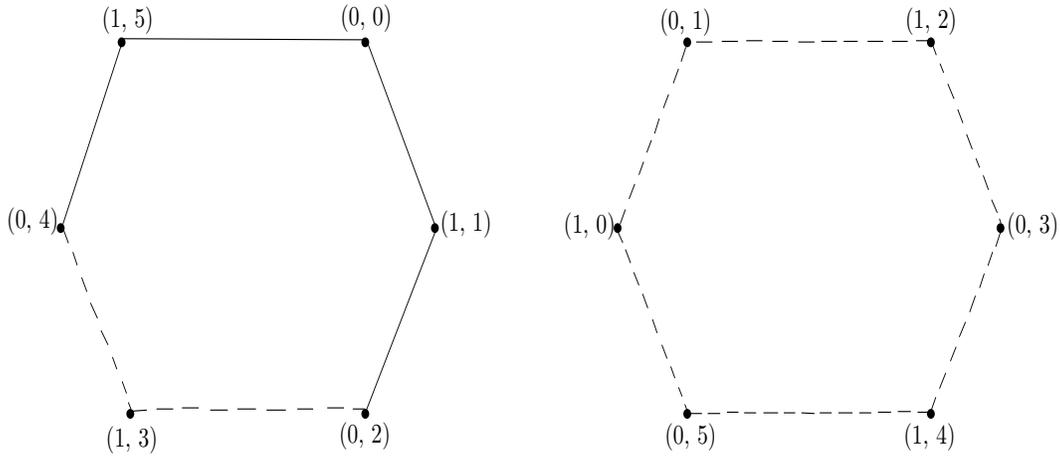} 
\caption{ Two connected components of $\Sigma$, left $\Sigma_{1}$, right $\Sigma_{2}$}\label{fig:1}
\end{center}
\end{figure}

\begin{figure}[h]
\begin{center}
\includegraphics[width=14cm , height=8cm]{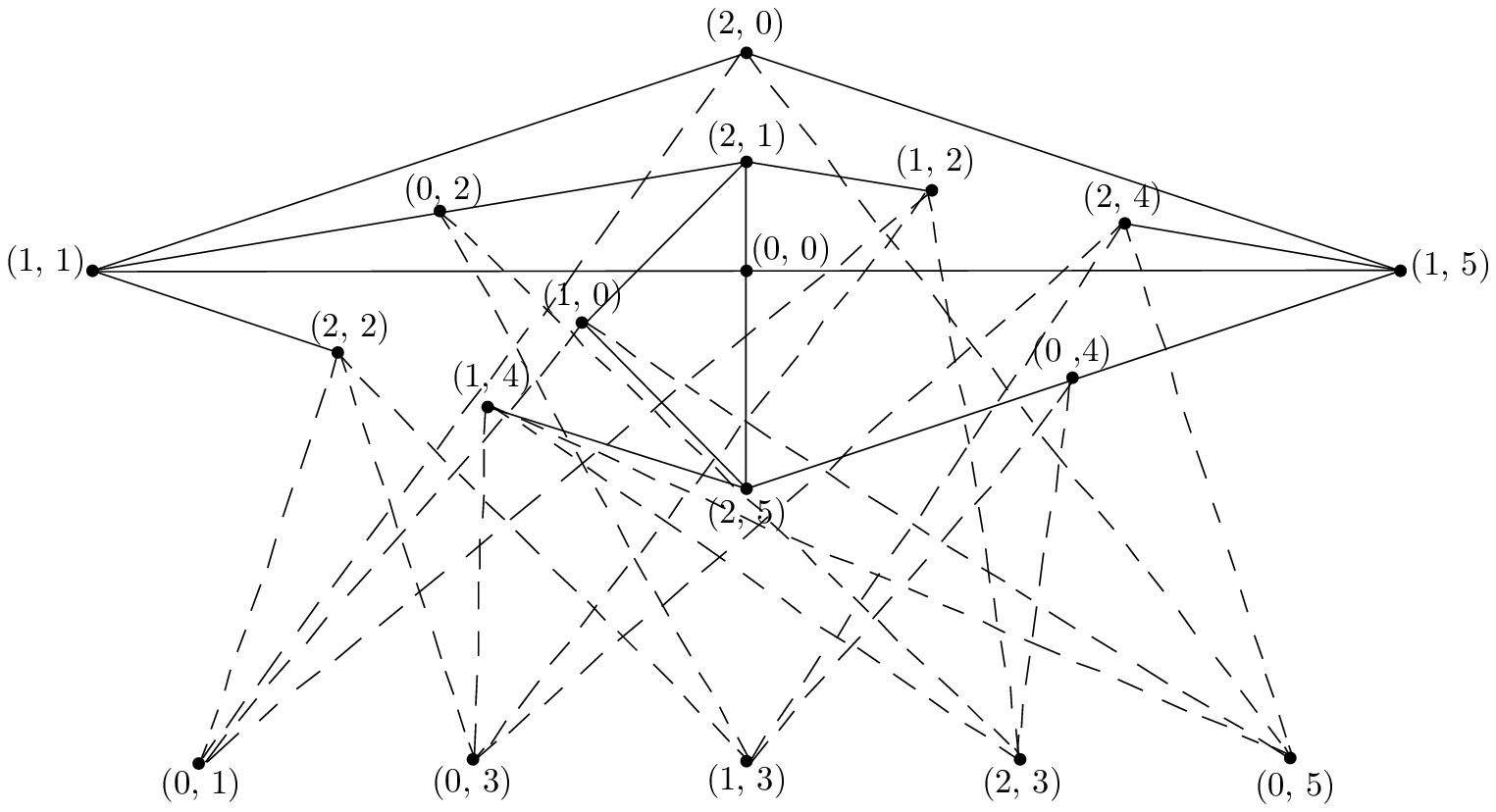} 
\caption{ sigraph $\Sigma^{'}$}\label{fig:2}
\end{center}
\end{figure}

When $\Gamma=\Cay(\mathbb{Z}_{2}\times \mathbb{Z}_{2}, \Phi)$, sigraph $\Sigma$ 
has two connected components. One of the components has a positive edge and the other 
component has a negative edge. (According to the explanation given in the last paragraph of the
 introduction, $V(\Sigma_{1})=\{(0, 0), (1, 1)\}$ and $V(\Sigma_{2})=\{(0, 1), (1, 0)\}$. Since 
 $(1, 1)\in \Phi$ and $(0, 0)$ is adjacent to $(1, 1)$, so, by definition of $\sigma$, the edge $(0, 0)(1, 1)$
  is positive. Also, since the vertices of $V(\Sigma_{2})$ are not in $\Phi$ and $(0, 1)$ is adjacent to $(1, 0)$,
   hence $(0, 1)(1, 0)$ is negative edge). Note that the tree has no cycle and using Theorem 1.6 \cite{s15}, 
we can conclude that each tree is balanced. Since $\Sigma$ is the union of two trees, so it is 
balanced. Also $L(\Sigma)$ has only two vertices so we may consider it is a balansed graph. Therefore
 in the proof of Theorems \ref{p3.7}, \ref{p3.9}, this case is not considered. 
\begin{prop}\label{p3.1}
Let $\Sigma=(\Gamma, \sigma)$, where 
$\Gamma=\Cay(\mathbb{Z}_{p_{1}}\times \mathbb{Z}_{p_{1}^{\alpha_{1}}p_{2}^{\alpha_{2}} \ldots p_{k}^{\alpha_{k}}}, \Phi)$
 and one of the prime factors is $2$. Then $\Sigma$ has $|\Phi|^2$ positive edges.
\begin{proof}
Since one of $p_{i}$ is $2$, so without loss of generality we may assume that $p_{1}=2$.
 In this case $\Sigma$ is a disconnected sigraph with exactly 
two connected components $\Sigma_{1}=(\Gamma_{1}, \sigma)$ and 
$\Sigma_{2}=(\Gamma_{2}, \sigma)$, where 
$V(\Sigma_{1})=\{(1, v)| v~is~odd\}\cup \{(0, v)| v~is~even\}$ and
$V(\Sigma_{2})=\{(0, v)| v~is~odd\}\cup \{(1, v)| v~is~even\}$. 
Note that $\Sigma$ is $|\Phi|$-regular hence $|E(\Sigma_{1})|=|E(\Sigma_{2})|=\dfrac{1}{4}.|\Phi|.p_{1}^{\alpha_{1}+1}p_{2}^{\alpha_{2}} \ldots p_{k}^{\alpha_{k}}$.

Suppose $(u, v)(u^{'}, v^{'})$ is an arbitrary edge of $\Sigma_{1}$, 
where $(u, v)\in \{(1, v)| v~is~odd\}$ and $(u^{'}, v^{'})\in \{(0, v)| v~is~even\}$.
 Since $\Phi\subseteq\{(1, v)| v~is~odd\}$ and by the definition of $\sigma$ the edge 
 that has at least one end in $\Phi$ is positive, so a number of edges in $\Sigma_{1}$
 are positive.  Also each of the members of $\Phi$ is adjacent to the $|\Phi|$ vertices of
  $\{(0, v)| v~is~even\}$. Hence we have $|\Phi|^{2}$ positive edges. Note that if 
  $\Gamma=\Cay(\mathbb{Z}_{2}\times \mathbb{Z}_{2^{\alpha_{1}}}, \Phi)$
   then $\Phi=\{(1, v)| v~is~odd\}$. Therefore $\Sigma_{1}$ is all-positive component of 
   sigraph $\Sigma$ (since the edges of $\Sigma_{1}$ have one end in $\{(1, v)| v~is~odd\}$
    and the other in $\{(0, v)| v~is~even\}$, so by definition of $\sigma$ are all positive). Also
     $|E(\Sigma_{1})|=\dfrac{1}{4}.|\Phi|.p_{1}^{\alpha_{1}+1}
   =|\Phi|.2^{\alpha_{1}-1}=|\Phi|^{2}$, since the number of odd numbers in $\{1, 2, 3, \ldots, 2^{\alpha_{1}}\}$ 
   is equal to $2^{\alpha_{1}-1}=|\Phi|$. 
   
 Now suppose that $(u, v)(u^{'}, v^{'})\in E(\Sigma_{2})$, where $(u, v)\in\{(0, v)| v~is~odd\}$ 
 and $(u^{'}, v^{'})\in \{(1, v)| v~is~even\}$. Clearly $u$ and $v^{'}$ are multiples of $2$.
  So  $(u, v)$ and $(u^{'}, v^{'})$ are not in $\Phi$. This implies that $\Sigma_{2}$ is the 
  all-negative component of $\Sigma$ (since the vertices of $\Sigma_{2}$ are not in $\Phi$, 
  so all the edges in $\Sigma_{2}$ are negative).    
  
   Next assume that $p_{1}\geq 3$. According to the definition of a Cayley graph, vertex $(0, 0)$ is adjacent to
   all vertices of $\Phi$. Also each vertex $\Phi$ is adjacent to $|\Phi|-1$ of other vertices of $\Sigma$.
    By definition of $\sigma$, all edges between these vertices are positive and their number is equal
     to $|\Phi|+|\Phi|(|\Phi|-1|)=|\Phi|^{2}$. We represent the set of all these vertices with $V_{1}$ such that 
  $V_{1}=\{(0, v)|v ~is~even\}\cup \{(u, v)| u\in \{1, 2, \ldots, p_{1}-1\}~and~v~is~not~odd$-$multiple~of~prime~factors\}$.
  Consider the vertices of $V_{2}$, where 
  
  $V_{2}=\{(0, v)|v ~is~odd\}\cup \{(u, v)| u\in \{1, 2, \ldots, p_{1}-1\}~and~v~is~odd$-$multiple~of~prime~factors\}$.
   Thus $V=V_{1}\cup V_{2}$. Assume that $(u, v)$ and $(u^{'}, v^{'})$ are 
   arbitrary vertices of $V_{2}$. Since $v$ and $v^{'}$ are both odd so $v-v^{'}$
    is even, hence $(u, v)-(u^{'}, v^{'})\notin \Phi$. Therefore none of the vertices 
    $V_{2}$ are adjacent to each other. Let's assume now $(u,v)\in V_{2}$ is adjacent
     to $(u^{'}, v^{'})\in V_{1}$. Because $v$ is an odd, so $v^{'}$ should be even. 
     Thus $(u^{'}, v^{'})\notin \Phi$. Moreover it is easy to see  that the vertices of
      $V_{2}$ are not in $\Phi$. This implies that all edges incident with the vertices of $V_{2}$ 
      are negative and their number is $|V_{2}|\times |\Phi|=\dfrac{n|\Phi|}{2}-|\Phi|^{2};~
       n=p_{1}^{\alpha_{1}+1}p_{2}^{\alpha_{2}} \ldots p_{k}^{\alpha_{k}}$.  
       \end{proof}
\end{prop}
  \begin{ex}\label{p3.2}
  Let $\Sigma=(\Cay(\mathbb{Z}_{2}\times \mathbb{Z}_{2\times 3}, \Phi), \sigma)$ and $\Sigma^{'}=(\Cay(\mathbb{Z}_{3}\times \mathbb{Z}_{2\times 3}, \Phi), \sigma)$ which are
   shown in Figures \ref{fig:1} and \ref{fig:2}, respectively. Positive edges with straight line and negative
    edges are shown with dash.
     \end{ex}

\begin{figure}[h]
\begin{center}
\includegraphics[width=7cm , height=6cm]{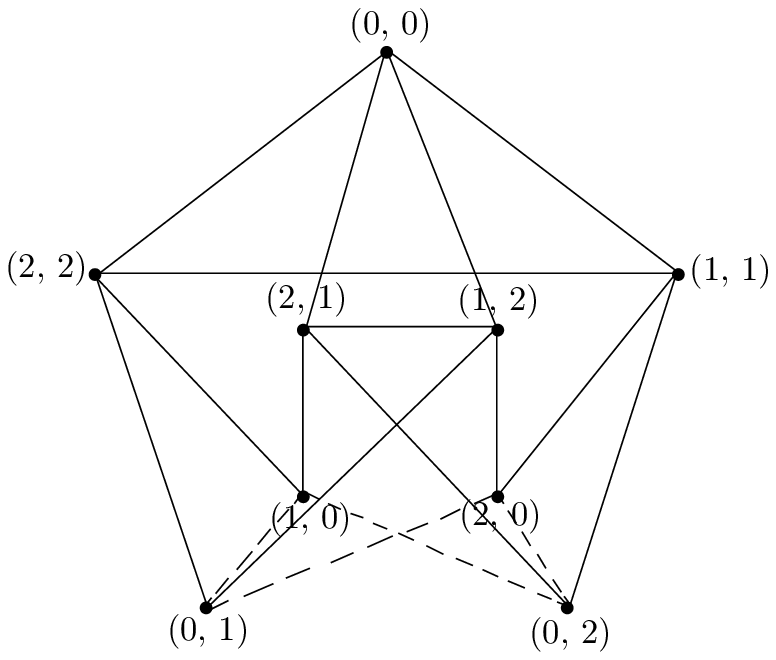} 
\caption{ sigraph $\Sigma$}\label{fig:3}
\end{center}
\end{figure}

  \begin{lem}\label{p3.3}
   Let $\Sigma=(\Gamma, \sigma)$, where 
$\Gamma=\Cay(\mathbb{Z}_{p}\times \mathbb{Z}_{p^{\alpha}}, \Phi)$, $\alpha\geq 1$ and $p\geq 3$ is prime number. 
Then the number of negative edges is $p^{\alpha-1}|\Phi|$.
\begin{proof}
Assume that $(u, v)(u^{'}, v^{'})$ is an arbitrary edge of $\Sigma$. 
By the definition of $\sigma$, we have negative edge only when $u, v^{'}$ or
 $u^{'}, v$ are multiples of $p$. Without loss of generality we suppose that $u, v^{'}$ 
are multiples of $p$. Since in $\mathbb{Z}_{p}$ only zero is multiple of $p$, hence
 $u=0, v\in \varphi_{p^{\alpha}}$ and also $u^{'}\in \varphi_{p}$ and $v^{'}$ is 
multiple of $p$ less than $p^{\alpha}$. We know that $\varphi_{p^{\alpha}}=p^{\alpha-1}(p-1)$ 
so $(u, v)$ can have $p^{\alpha-1}(p-1)$ cases. Moreover $u^{'}\in\{1, 2, \ldots, p-1\}$ and we
 can consider $p^{\alpha-1}$ cases for $v^{'}$, because the number of multiples of $p$ less than
  $p^{\alpha}$ is equal to $p^{\alpha-1}$. Therefore $(u^{'}, v^{'})$ can have $p^{\alpha-1}(p-1)$ cases.
  This implies that $\Sigma$ has 
  $(p^{\alpha-1}(p-1))^{2}=p^{\alpha-1}.|\varphi_{p}|.|\varphi_{p^{\alpha}}|=p^{\alpha-1}|\Phi|$ negative edges.
\end{proof}
\end{lem}
\begin{ex}\label{p3.4}
Let $\Sigma=( \Cay(\mathbb{Z}_{3}\times \mathbb{Z}_{3}, \Phi), \sigma)$. This sigraph, which is shown in figure \ref{fig:3}, has $4$ negative edges. 
\end{ex}
  \begin{remark}\label{p3.5}
  Suppose $\lambda$ is the length of the longest sequence
of consecutive integers in $\mathbb{Z}_{n}$, each of which shares a prime factor with 
$n$. Assume that $n=pq$, where $p, q$ are prime numbers and $p, q\geq 3$. In 
$\mathbb{Z}_{pq}$ we have $\lambda=2$. Now we prove that $\lambda$ appears exactly in two places. 

Let $kp$ and $k'q$ be the first consecutive integers in $\mathbb{Z}_{pq}$, which are
 multiples of $p$ and $q$, respectively. It is clear that
$pq-k'q$ and $pq-kp$ are the final consecutive
integers in $\mathbb{Z}_{pq}$. Suppose that $fp$ and $tq$ are other consecutive integers in $\mathbb{Z}_{pq}$,
where $f>k$ and $t>k'$. Then $(fp-tq)-(k'q-kp)=0$ and therefore $(f+k)p=(t+k')q$.
Hence there exists an integer $r$ with $f+k=rq$ and $t+k'=rp$. Since $1\leq k<f\leq q-1$
and $1\leq k'<t\leq p-1$, we have $r=1$. So $f=q-k$ and $t=p-k'$. Therefore $fp$ and $tq$
are the final consecutive integers in $\mathbb{Z}_{pq}$.
  \end{remark}
  
  \begin{lem}\label{p3.6}
Let $\Sigma=(\Gamma, \sigma)$, where 
$\Gamma=\Cay(\mathbb{Z}_{p}\times \mathbb{Z}_{pq}, \Phi)$ and $p, q\geq 3$. 
Then $\Sigma$ has $|\Phi|(2p+q-3)$ negative edges.
\begin{proof}
Let's consider the arbitrary edge $(u, v)(u^{'}, v^{'})$ of $\Sigma$. Only when we have a 
negative edge that the vertices $(u, v)$ and $(u^{'}, v^{'})$ are not in $\Phi$. Also
 this is possible if at each of the vertices $(u, v), (u^{'}, v^{'})$ at least one of theire elements
  multiples of one of the prime factors $p$ or $q$ .The following cases are the only possible
   cases for the formation of negative edges. Note that we show nonzero multiples of $p$ and $q$ by 
   $\alpha p$ and $\beta q$, respectively.  
   
   \textbf{Case(i)} $(u, v)\in \{(0, v)| v \in \varphi_{pq}\}$, 
 $(u^{'}, v^{'})\in\{(u^{'}, v^{'})|u^{'}\in \{1, 2, \ldots, p-1\}, v^{'}=0~ or ~v^{'}=
 \alpha p, v-v^{'}\neq \beta q~or~v^{'}=\beta q, v-v^{'}\neq \alpha p\}$.
 
 Clearly there are $|\varphi_{pq}|$ cases for $(u, v)$. For $(u^{'}, v^{'})$ if $v^{'}=0$ 
 then we have $p-1$ cases and if $v^{'}=\alpha p, v-v^{'}\neq \beta q$ then by Remark \ref{p3.5}, we have $(p-1)(q-2)$
  cases and also if $v^{'}=\beta q, v-v^{'}\neq \alpha p$ then we have $(p-1)(p-2)$ cases for $(u^{'}, v^{'})$.
   Therefore in this case , the number of negative edges is $|\varphi_{pq}|((p-1)+(p-1)(q-2)+(p-1)(p-2))=|\Phi|(p+q-3)$.
   
   \textbf{Case(ii)} $(u, v)\in \{(0, v)| v=\alpha p\}$, 
 $(u^{'}, v^{'})\in\{(u^{'}, v^{'})|u^{'}\in \{1, 2, \ldots, p-1\}, v^{'}=\beta q\}$.
 
 In this case we have $q-1$ cases for $(u, v)$ and  $(p-1)(p-1)$ cases for $(u^{'}, v^{'})$. Hence
  there are $(q-1)(p-1)(p-1)=|\Phi|$ cases for $(u, v)(u^{'}, v^{'})$.
  
  \textbf{Case(iii)} $(u, v)\in \{(0, v)| v=\beta q\}$, 
 $(u^{'}, v^{'})\in\{(u^{'}, v^{'})|u^{'}\in \{1, 2, \ldots, p-1\}, v^{'}=\alpha p\}$. 
 
 Similar to Case(ii), 
 there are $(p-1)(p-1)(q-1)=|\Phi|$ cases for $(u, v)(u^{'}, v^{'})$.
 
 \textbf{Case(iv)} $(u, v)\in \{(u, v)| u\in \{1, 2, \ldots, p-1\}, v=\alpha p\}$, 
 $(u^{'}, v^{'})\in\{(u^{'}, v^{'})|u^{'}\in \{1, 2, \ldots, p-1\}, u^{'}\neq u, v^{'}=\beta q\}$.  
 There are $p-1,~q-1,~p-2$ and $p-1$ cases for $u,~v,~u^{'}$ and $v^{'}$, respectively. Therefore 
 we have $(p-1)(q-1)(p-2)(p-1)=|\Phi|(p-2)$ cases for $(u, v)(u^{'}, v^{'})$. 
 
 By considering all the above cases, the number of negative edges in $\Sigma$ is equal to:
 
  $|\Phi|(p+q-3)+|\Phi|+|\Phi|+|\Phi|(p-2)=|\Phi|(2p+q-3)$.
\end{proof}
\end{lem}

  \begin{thm}\label{p3.7}
  Let $\Sigma=(\Gamma, \sigma)$, where 
$\Gamma=\Cay(\mathbb{Z}_{p_{1}}\times 
\mathbb{Z}_{p_{1}^{\alpha_{1}}p_{2}^{\alpha_{2}} \ldots p_{k}^{\alpha_{k}}}, \Phi)$.
 Then $\Sigma$ is balanced if and only if one of the prime factors is $2$.
\begin{proof}
\textbf{Necessity}: Let $\Sigma$ is balanced. On contrary, assume that $p_{i}\geq 3$
 for any $ i=1, 2, \ldots, k$. Sinse $p_{i}\geq 3$ hence both numbers $1$ and $2$ 
belong to $\varphi_{p_{1}}$ and $\varphi_{p_{1}^{\alpha_{1}}p_{2}^{\alpha_{2}} \ldots p_{k}^{\alpha_{k}}}$. 
This implies that $\Sigma$ contains the cycle $C=((0, 1), (1, 2), (2, 0), (0, 1))$. 
Note that $(1, 2)\in\Phi$ hence $(0, 1)(1, 2)$ and $(1, 2)(2, 0)$ are positive edges. 
Moreover $(0, 1)$ and $(2, 0)$ are not in $\Phi$, so $(2, 0)(0, 1)$ is negative edge. 
 We conclude that $C$ is a negative cycle. Therefore, by definition of balanced sigraph, $\Sigma$ is 
 unbalanced. Which is in contradiction with the assumption.
 
 \textbf{Sufficiency}: Suppose one of the prime factors is $2$. Assume first that $p_{1}=2$.
  Then by proof of Proposition \ref{p3.1}, $\Sigma$ has exactly two connected components $\Sigma_{1}$
   and $\Sigma_{2}$ where $\Sigma_{2}$ is all-negative. We show that every cycle in $\Sigma_{2}$
    is positive. Suppose $C=(x_{1}, x_{2}, \ldots, x_{m}, x_{1})$ is an arbitrary cycle in $\Sigma_{2}$.
     Without loss of generality, let $x_{1}\in \{(0, v)| v~is~odd\}$. Since neither of the vertices of 
     $\{(0, v)| v~is~odd\}$ is not adjacent to each other, so $x_{2}, x_{m}\in \{(1, v)| v~is~even\}$.
      By continuing this process, it is easy to see that $C$ contains an even number of negative edges.
       Which implies that $\Sigma_{2}$ is balanced.
 
 Now we prove that $\Sigma_{1}$ is also balanced. Let 
 $\Gamma=\Cay(\mathbb{Z}_{2}\times\mathbb{Z}_{2^{\alpha}}, \Phi)$. 
 Then $\Phi$ is consists of all the members of $\{(1, v)| v~is~odd\}$ and since 
 $\{(1, v)| v~is~odd\}\subseteq V(\Sigma_{1})$ Thus $\Sigma_{1}$ is all-positive.
  Clearly in this case $\Sigma_{1}$ is balanced. Assume that $\Gamma=\Cay(\mathbb{Z}_{2}\times 
\mathbb{Z}_{2^{\alpha_{1}}p_{2}^{\alpha_{2}} \ldots p_{k}^{\alpha_{k}}}, \Phi)$. 
Consider the arbitrary cycle $C^{'}=(x^{'}_{1},x^{'}_{2},\ldots, x^{'}_{n}, x^{'}_{1} )$. 
The edges in $\Sigma_{1}$ are such that each edge has a vertex in $\{(0, v)| v~is~even\}$ 
and a vertex in $\{(1, v)| v~is~odd\}$. Now without loss of generality assume that
 $x^{'}_{1}\in\{(0, v)| v~is~even\}$ hence $x^{'}_{2}, \in \{(1, v)| v~is~odd\}$.
 Hence we have two following cases. If $x^{'}_{2}\in\Phi$ then all edges joint to $x^{'}_{2}$ 
 will be positive. So $x^{'}_{1}x^{'}_{2}$ and $x^{'}_{2}x^{'}_{3}$ are positive.
  But if $x^{'}_{2}\notin\Phi$, because $x^{'}_{1}\notin\Phi$, then $x^{'}_{1}x^{'}_{2}$ 
  is negative. Since $x^{'}_{3}\in\{(0, v)| v~is~even\}$ hence $x^{'}_{2}x^{'}_{3}$ is also 
  negative. This implies that after each negative edge in $C^{'}$ we will have another negative edge. 
  Therefore, the number of negative edges in $C^{'}$ is even. Thus $\Sigma_{1}$ is balanced. 
  
  Next consider the case where $p_{1}\geq 3$. In this case $\Sigma$ is connected sigraph. 
  Suppose $C^{''}$ is an arbitrary cycle in $\Sigma$. If all the edges in $C^{''}$ are positive
   then the cycle is positive. Now assume that $C^{''}$ contains a negative edge $(u, v)(u^{'}, v^{'})$. 
   Using the proof of Proposition \ref{p3.1} and Figure \ref{fig:2} of Example \ref{p3.2}, $(u, v)(u^{'}, v^{'})$ has a vertex in
    $V_{1}$ and other vertex in $V_{2}$. Without loss of generality suppose that $(u, v)\in V_{1}$ and
    $ (u^{'}, v^{'})\in V_{2}$. Because neither of the vertices of $V_{2}$ is adjacent to each other, so 
    $(u^{'}, v^{'})$ must be connected to one of the vertices of $V_{1}$, for example $(u^{''}, v^{''})$. 
    Since all the edges that have a vertex in $V_{2}$ are negative, thus $(u^{'}, v^{'})(u^{''}, v^{''})$ is negative edge. 
    This implies that, in order to form cycle $C^{''}$, for each negative edge in $C^{''}$, another negative 
    edge is required. We conclude that the number of negative edges in $C^{''}$ is even. Therefore $\Sigma$ is balanced.  
\end{proof}
  \end{thm}
  \begin{thm}\cite{s1}\label{p3.8}
  For a sigraph $S$, its line sigraph $L(S)$ is balanced if and only if the
following conditions hold:
\item[1)] for any cycle $Z$ in $S$,\\
(a) if $Z$ is all-negative, then $Z$ has even length;\\
(b) if $Z$ is heterogeneous, then $Z$ has even number of negative sections with
even length;
\item[2)] for $v\in S$, if $d(v) > 2$, then there is at most one negative edge incident at $v$ in
$S$.
\end{thm}
  \begin{thm}\label{p3.9}
  Let $\Sigma=(\Gamma, \sigma)$, where 
$\Gamma=\Cay(\mathbb{Z}_{p_{1}}\times 
\mathbb{Z}_{p_{1}^{\alpha_{1}}p_{2}^{\alpha_{2}} \ldots p_{k}^{\alpha_{k}}}, \Phi)$.
 Then $L(\Sigma)$ is balanced if and only if one of the prime factors is $2$.
  \begin{proof}
  \textbf{Necessity}: Let $L(\Sigma)$ is balanced. Assume that the conclusion is false. Suppose $p_{i}\geq 3$
 for any $ i=1, 2, \ldots, k$. It is easy to see that vertex $(0, 1)$ is adjacent to vertices 
  $(1, 0), (2, 0), \ldots, (p_{1}-1, 0)$. By definition $\Phi$, the vertices $(0, 1), (1, 0), (2, 0), \ldots, (p_{1}-1, 0)$ are not in
   $\Phi$. Hence $d^{-}((0, 1))\geq 2$. By the Theorem \ref{p3.8}, condition (ii) does not hold. Therefore
   $L(\Sigma)$ is unbalansed. Which is in contradiction with the assumption. Thus one of the prime factors is $2$.
  
   \textbf{Sufficiency}: Suppose one of the prime factors is $2$. By Theorem \ref{p3.7}, $\Sigma$ is 
 balanced so, according to Lemma 1.5\cite{s15} and Theorem 1.6\cite{s15}, $\Sigma$ is switching 
 equivalent to a graph whose all edges are positive. Hence $\Sigma$ can be considered all-positive. 
 Therefore $L(\Sigma)$ is balanced.
  \end{proof}
  \end{thm}
  \section{Clusterability and sign-compatibility of $\Sigma$}\label{sec:3}
\begin{figure}[h]
\begin{center}
\includegraphics[width=6cm , height=2cm]{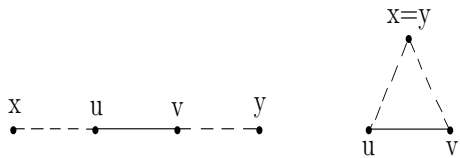} 
\caption{ Forbidden subsigraphs for a sign-compatible sigraph.}\label{fig:4}
\end{center}
\end{figure}

  In this chapter two properties namely clusterability and sign-compatibility for $\Sigma$ are investigated. 
  \begin{thm}\cite{s2}\label{p3.10}
  Let $S$ be any signed graph. Then $S$ has a clustering if and only if $S$ contains 
  no cycle having exactly one negative line.
  \end{thm}
   \begin{thm}\cite{s3}\label{p3.12}
  A sigraph $S$ is sign-compatible if and only if $S$ does not contain a subsigraph
that is isomorphic to either of the following two sigraphs formed from the path graph
$P4 = (x, u, v, y)$; $S1$ with both edges $xu$ and $vy$ negative and the edge $uv$ positive and $S2$
obtained from $S1$ by identifying the vertices $x$ and $y$ (see Figure \ref{fig:4} ).
   \end{thm}
  \begin{thm}\label{p3.11}
  Let $\Sigma=(\Gamma, \sigma)$, where 
$\Gamma=\Cay(\mathbb{Z}_{p_{1}}\times 
\mathbb{Z}_{p_{1}^{\alpha_{1}}p_{2}^{\alpha_{2}} \ldots p_{k}^{\alpha_{k}}}, \Phi)$.
 Then $\Sigma$ is clusterable if and only if $\Sigma$ is balanced. 
  \begin{proof}
  \textbf{Nesssity}: Suppose $\Sigma$ is clusterable. Assume that the conclusion is false. 
  Suppose $\Sigma$ is unbalanced. Using Theorem \ref{p3.7}, $p_{i}\geq 3$ for any 
  $ i=1, 2, \ldots, k$. Because $p_{i}\geq 3$ so $(2, 2)\in \Phi$ and also it is adjacent to
   vertices $(0, 1)$ and $(1, 0)$. We now consider the cycle $C=((0, 1), (2, 2), (1, 0), (0, 1))$ 
   in $\Sigma$. By the definition of $\sigma$, we have $\sigma((0, 1)(2, 2))=\sigma((2, 2)(1, 0))=+$ 
   and $\sigma((1, 0)(0, 1))=-$. Hence $C$ is a cycle with exactly one negative edge. Therefore, according 
   to Theorem \ref{p3.10}, $\Sigma$ can not be clusterable, a contradiction to the hypothesis.
   Hence one of the prime factors is $2$ so $\Sigma$ is balanced.
   
   \textbf{Sufficiency}: Suppose $\Sigma$ is balanced. Thus all cycles of $\Sigma$ are positive. 
   So the number of negative edges in them is even. Hence it can not include the cycle with a single 
   negative edge. Therefore, according to Theorem \ref{p3.10}, $\Sigma$ is clusterable. 
  \end{proof}
  \end{thm}
 
  \begin{thm}
   Let $\Sigma=(\Gamma, \sigma)$, where 
$\Gamma=\Cay(\mathbb{Z}_{p_{1}}\times 
\mathbb{Z}_{p_{1}^{\alpha_{1}}p_{2}^{\alpha_{2}} \ldots p_{k}^{\alpha_{k}}}, \Phi)$ and $p_{i}\geq 2$ for any 
  $ i=1, 2, \ldots, k$. Then $\Sigma$ is sign-compatible.
  \begin{proof}
  Assume $\Sigma$ is sign-incompatible. Then by Theorem \ref{p3.12}, $\Sigma$ contains a subsigraph 
  isomorphic to either $S1$ or $S2$. Let $P=((u, v), (u^{'}, v^{'}), (u^{''}, v^{''}), (u^{'''}, v^{'''}))$ is 
  a path in $\Sigma$ such that $\sigma((u, v)(u^{'}, v^{'}))=\sigma((u^{''}, v^{''})(u^{'''}, v^{'''}))=-$ and
   $\sigma((u^{'}, v^{'})(u^{''}, v^{''}))=+$. According to the definition of $\sigma$, since
   $(u^{'}, v^{'})(u^{''}, v^{''})$ is positive edge hence two cases may occure. 
  
  If $(u^{'}, v^{'})$ and $(u^{''}, v^{''})$ belong to $\Phi$ then all the edges on $P$ are positive, 
  a contradiction to hypothesis. Now if one of the vertices $(u^{'}, v^{'})$ and $(u^{''}, v^{''})$ 
  belong to $\Phi$. Without loss of generality assume that $(u^{'}, v^{'})\in \Phi$ and $(u^{''}, v^{''})\notin \Phi$.
   Then we have $\sigma((u, v)(u^{'}, v^{'}))=((u^{'}, v^{'})(u^{''}, v^{''}))=+$ and 
   $\sigma((u^{''}, v^{''})(u^{'''}, v^{'''}))=-$. Which is contradiction the hyposesis. Thus 
   $\Sigma$ does not contain a subsigraph isomorfic to $S1$.
  
  Now let $\Sigma$ contains a subsigraph isomorphic to $S2$. As proved for $S1$,
   according to the definition of $\sigma$, no positive edge can be placed between two 
   negative edges. Hence $\Sigma$ does not contain a subsigraph isomorfic to $S2$. 
   Therefore $\Sigma$ is sign-compatible. 
  \end{proof}
  \end{thm}




\end{document}